\documentclass[final,3p,times]{elsarticle}

\usepackage{mathrsfs}
\usepackage{amsfonts}
\usepackage{color}
\usepackage{mathrsfs}
\usepackage{mathtools}
\usepackage{stmaryrd}
\usepackage{amsmath}
\usepackage{amssymb}
\usepackage{bbm}
\usepackage{mathrsfs}
\usepackage{graphicx}
\usepackage{enumerate}
\usepackage{amsthm}
\newcommand{\ba}{\begin{eqnarray}}
\newcommand{\ea}{\end{eqnarray}}
\newtheorem{thm}{Theorem}[section]

\newtheorem{lem}{Lemma}[section]
\newtheorem{defi}{Definition}[section]

\usepackage{amssymb}
\def\RR{\mathbb{R}}
\def\PP{\mathbb{P}}
\def\FF{\mathbb{F}}
\def\EE{\mathbb{E}}

\def\si{{\sigma}}

\def\De{{\Delta}}

\def\Om{{\Omega}}
\def\om{{\omega}}
\def\al{{\alpha}}

\def\De{{\Delta}}

\def\si{{\sigma}}

\def\EE{\mathbb{ E}}

\def\si{{\sigma}}
\def\al{{\alpha}}

\renewcommand{\d}{d}
\begin{document}

\begin{frontmatter}



\title{\Large \bf The stochastic nonlinear Schr\"{o}dinger equations driven by pure jump noise}


\author[inst1]{Jian Wang}
\ead{wg1995@mail.ustc.edu.cn}
\affiliation[inst1]{organization={School of Mathematical Sciences},%Department and Organization
            addressline={University of Science and Technology of China}, 
            city={Hefei},
            postcode={230026}, 
            %state={State One},
            country={China}}

\author[inst1]{Jianliang Zhai}
\ead{zhaijl@ustc.edu.cn}
\author[inst2]{Jiahui Zhu}
\ead{jiahuizhu@zjut.edu.cn}
\affiliation[inst2]{organization={School of Science},%Department and Organization
            addressline={Zhejiang University of Technology}, 
            city={Hangzhou},
            postcode={310019}, 
         %   state={State Two},
            country={China}}

\begin{abstract}
In this paper, we establish the existence and uniqueness of solutions  of stochastic nonlinear Schr\"{o}dinger equations with additive jump noise in $L^2(\RR^d)$. Our results cover all either focusing or defocusing nonlinearity in the full subcritical range of exponents as in the deterministic case.
\end{abstract}




\end{frontmatter}

 
\section{Introduction and motivation}
In this paper, we consider the following stochastic nonlinear Schr\"{o}dinger equations driven by pure jump noise:
\begin{eqnarray}\label{NLS}
dX(t) = i[ \De X(t) - \lambda |X(t)|^{\al -1}X(t)]dt + dL(t), \qquad X(0) = x \in L^2(\RR^d),
%\int_{Z_1^c} z \tilde{N}(dz, dt) +\int_{Z_1} z N(dz, dt), \notag\\
\end{eqnarray}
where $\lambda \in \{-1, 1\}$ and $1 < \al< 1+\frac{4}{d}$ and $L=(L(t))_{t\geq 0}$ is an $L^2(\RR^d)$-valued pure jump L\'evy process defined as
$$L(t) = \int_0^t \int_B z \tilde{N} (dz,dt),$$ where $B = \{z \in L^2(\RR^d): 0 < \|z\|_{L^2(\RR^d)} \leq 1\}$. Here, $N$ represents a time homogeneous Poisson random measure over $(\mathcal{B}(L^2(\RR^d)), \mathcal{B}(\mathbb{R}^{+}))$ with $\si$-finite intensity measure  $\nu$ satisfying $\int_{L^2(\RR^d) \backslash \{0\}}\|z\|_{L^2(\RR^d)}^2\wedge1\nu(dz)<\infty$ and $\tilde{N} (dz,dt) = N(dz,dt) - \nu(dz)dt$ denotes the corresponding compensated Poisson random measure.

The nonlinear Schr\"{o}dinger equation appears in a variety of applications in nonlinear physics, including nonlinear optics, nonlinear water propagation, nonlinear acoustics, quantum condensates etc. In recent years, existence and uniqueness of solutions for the stochastic nonlinear Schr\"{o}dinger equation with Gaussian noise have been investigated by many authors \cite{BRZ-14,BRZ-17,deB+Deb-99, deB+Deb-03,BHW-19,Hornung}. When it comes to L\'evy noise,  the current literature  is still very limited.  In \cite{BHW-19}, Brze\'zniak et al. proved the existence of a martingale solution of the stochastic nonlinear Schr\"{o}dinger equation with a multiplicative jump noise in the Marcus canonical form in $H^1$ on compact manifolds by using the Gaedo-Galerkin method. De Bouard and Hausenblas  investigated the existence of martingale solutions of the nonlinear Schr\"{o}dinger equation with a L\'evy noise with infinite activity in \cite{deB+Hau-19} and pathwise uniqueness was studied in a separate paper \cite{deB+Hau+Ond-19} with Ondrejat. Recently, Brze\'zniak et al. established a new version of the stochastic Strichartz estimate for the stochastic convolution driven by a jump noise in \cite{BLZ}. By applying the stochastic Strichartz estimates in a fixed point argument, they proved the existence and uniqueness of a global solution to stochastic nonlinear Schr\"{o}dinger equation with a Marcus-type jump noise in $L^2(\mathbb{R}^d)$ with either focusing or defocusing nonlinearity in the full subcritical range of exponents.

%The authors \cite{BLZ} prove the existence and uniqueness of a global solution to stochastic nonlinear Schr\"{o}dinger equation in $L^2(\RR^d)$ with either focusing or defocusing nonlinearity in the full subcritical range of exponents as in the determ
The purpose of this paper is to prove the existence and uniqueness of mild solutions for the stochastic nonlinear Schr\"{o}dinger equation \eqref{NLS} with the additive jump noise.  By means of the deterministic and stochastic Strichartz's estimates due to Brze\'zniak et al. from \cite{BLZ}, we apply the classical truncation procedure of the nonlinearities and use the well-known fixed point argument to construct a local solution. One main ingredient of establishing global solution in \cite{BLZ} is the mass conservation for the stochastic nonlinear Schr\"{o}dinger equaiton. Unlike the $L^2(\RR^d)$-norm conservation of solutions in the literature \cite{BLZ}, we establish some uniform $L^2(\RR^d)$-norm estimate of the local solution and selects some specific stopping times that we believe are critical to the proof. Let us formulate our main result of this paper.

\begin{thm}\label{main}
  Let $p\geq 2$, $1<\al<1+\frac{4}{d}$, $r=\al +1$  such that  $(p,r)$ is an admissible pair. For any $x \in  L^2(\RR^d)$,
there exists a unique global mild solution  $X^x=(X^x(t),t\in[0,\infty))$ of (\ref{NLS}) such that
\[
X^x\in D([0,\infty);L^2(\RR^d))\cap L^p(0,T;L^r(\RR^d)),\ \mathbb{P}\text{-a.s.}\ \omega\in\Omega.
\]
          Moreover, we have for all $t\in[0,T]$ and any $q \geq 2$
\begin{equation}\label{L2 norm}
               \EE \sup_{t \in [0,T]}\|X(t)\|_{L^2(\RR^d)}^q \leq C\big(q,T, \|x\|_{L^2(\RR^d)}\big). 
\end{equation}           
\end{thm}

\section{Setting and Strichartz estimates}
In this section, we introduce some notations, assumptions and solution concepts and
recall deterministic and stochastic Strichartz estimates, which will be used to prove the well-posedness of the solution of the truncated equation.

For $t>0$, let us denote $D(0,t;L^2(\RR^d))$ the space of all right continuous functions with left-hand limits from $[0,t]$ to $L^2(\RR^d)$ and
\begin{align}
Y_{t}:=L^{\infty}(0,t;L^2(\RR^d))\cap L^p(0,t;L^r(\RR^d)).
\end{align}
Then $Y_{t}$ is a Banach space with norm defined by
  \begin{align}
   \|u\|_{Y_{t}}:=\sup_{s\in[0,t]}\|u(s)\|_{L^2(\RR^d)}+\Big(\int_0^t\|u(s)\|_{L^r(\RR^d)}^p\d s\Big)^{\frac1p}, \quad  \text{for} \ u \in Y_t.
  \end{align}
  For the details we refer to (2.10) and (2.11) in \cite{BLZ}.
  
 Let $(\Omega,\mathcal{F}, \FF,\mathbb{P})$ be a filtered probability space satisfying the usual conditions, where $\FF = \{{\mathcal{F}}_t\}_{t\geq 0}$ is the filtration.  Let $\tau>0$ be a stopping time. We call $\tau$ an accessible stopping time if there exists an increasing sequence $(\tau_n)_{n\in\mathbb{N}}$ of stopping times such that $\tau_{n}<\tau$ and $\tau_{n} \nearrow \tau$ $\mathbb{P}$-a.s. as $n\rightarrow\infty$ and we call $(\tau_n)_{n\in\mathbb{N}}$ an approximating sequence for $\tau$.

   Let $M^p_{\mathbb{F}}(Y_{\tau}):=L^{p}(\Omega; L^{\infty}(0,\tau;L^2(\RR^n))\cap L^p(0,\tau;L^r(\RR^n)))$ denote the space of all $L^{2}(\RR^{n})\cap L^{r}(\RR^{n})$-valued $\mathbb{F}$-progressively measurable processes $u:[0,T]\times \Omega \rightarrow L^2(\RR^n)\cap L^{r}(\RR^{n})$  satisfying
  \begin{align*}
  \|u\|^p_{M^p_{\mathbb{F}}(Y_{\tau})}:=\EE\|u\|^p_{Y_\tau}=\EE\Big(\sup_{s\in[0,\tau]}\|u(s)\|^{p}_{L^2(\mathbb{R}^d)}+\int_0^{\tau}\|u(s)\|_{L^r(\mathbb{R}^d)}^p\d s\Big)<\infty.
  \end{align*}

Now we introduce the definitions of local solutions and maximal local solutions, see e.g. \cite{ZB} for more details.
  \begin{defi}
        A local mild solution to equation  (\ref{NLS}) is an $\FF$-progressively measurable process $X(t)$, $t\in[0,\tau)$, where $\tau$ is an accessible stopping time with an approximating sequence $(\tau_n)_{n\in\mathbb{N}}$ of stopping times such that for every $n\in\mathbb{N}$,
       \begin{enumerate}
        \item[(i)] $(X(t))_{t\in[0,\tau_{n}]}\in D(0,\tau_{n};L^{2}(\mathbb{R}^{d}))$, $\mathbb{P}$-a.s.;
\item[(ii)] $(X(t))_{t\in[0,\tau_{n}]}$ belongs to $M^p_{\mathbb{F}}(Y_{\tau_{n}})$;
\item[(iii)] for every $t\in [0,T]$, the following equality holds
\begin{align*}
          X(t\wedge \tau_{n})=&S_{t\wedge\tau_{n}}x -i\lambda\int_0^{t\wedge\tau_{n}}S_{t\wedge\tau_{n}-s}( |X(s)|^{\al-1}X(s))\,\d s +\int_0^{t\wedge\tau_{n}}\int_B 1_{[0, \tau_n]}(s)S_{t\wedge \tau_n-s}z \tilde{N} (dz, ds), \;\;\mathbb{P}\text{-a.s. }
        \end{align*}
              \end{enumerate}
 \end{defi}
   \begin{defi}
                 A local mild solution $X=(X(t))_{0\leq t<\tau}$ is called a maximal local mild solution if for any other local mild
   solution $(Y(t))_{t\in[0,\sigma) }$
   satisfying $\sigma\geq \tau$ a.s. and $Y|_{[0,\sigma)}$ is equivalent to $
  X$, one has  $\sigma=\tau$ a.s.. 
  
           A local mild solution $(X(t))_{t\in[0,\tau)}$ is a global mild solution if $\tau=T$, $\mathbb{P}$-a.s. and $u\in M^p_{\mathbb{F}}(Y_T)$.
 \end{defi}
% Let $(\Omega,\mathcal{F}, \FF,\mathbb{P})$ be a filtered probability space satisfying the usual conditions, where $\FF = \{{\mathcal{F}}_t\}_{t\geq 0}$ is the filtration. Let $L=(L(t))_{t\geq 0}$ be an $L^2(\RR^d)$-valued pure jump L\'evy process defined as
%$$L(t) = \int_0^t \int_B z \tilde{N} (dz,dt),$$ where $B = \{z \in L^2(\RR^d): 0 < \|z\|_{L^2(\RR^d)} \leq 1\}$. Here, $N$ represents a time homogeneous Poisson random measure over $(\mathcal{B}(L^2(\RR^d)), \mathcal{B}(\mathbb{R}^{+}))$ with $\si$-finite intensity measure  $\nu$ satisfying $\int_{L^2(\RR^d) \backslash \{0\}}\|z\|_{L^2(\RR^d)}^2\wedge1\nu(dz)<\infty$ and $\tilde{N} (dz,dt) = N(dz,dt) - \nu(dz)dt$ denotes the corresponding compensated Poisson random measure.

%Let $N: \cB(L^2(\RR^d)\times\RR^+) \times \Omega\rightarrow \bar{\NN} = \NN \cup \{0, \infty\}$ be the time homogeneous Poisson random measure with $\si$-finite intensity measure  $\nu$ satisfying $\int_{L^2(\RR^d) \backslash \{0\}}\|z\|_{L^2(\RR^d)}^2\wedge1\nu(dz)<\infty$. For the existence of such Poisson random measure, we refer the reader to \cite{IW 1989}. Again $\tilde{N} (dz,dt) = N(dz,dt) - \nu(dz)dt$ denotes the compensated Poisson random measure.
%Let $L(t)$ be an $L^2(\RR^d)$-valued pure jump L\'evy process with L\'evy measure $\nu$, i.e. $L(t) = \int_0^t \int_B z \tilde{N} (dz,dt)$, where $B = \{z \in L^2(\RR^d): 0 < \|z\|_{L^2(\RR^d)} \leq 1\}$.

Throughout the paper, the symbol $C$ will denote a positive generic constant whose value may change from place to place. If a constant depends on some variable parameters, we will put them in subscripts.

We now state the following deterministic and stochastic Strichartz estimates, we refer the reader to \cite{TC} and \cite[Propositions 2.2 and 2.6]{BLZ}.  Let $(S_t)_{t \in \RR}$ denote the group of isometries on $L^2(\RR^d)$ generated by $i \De$. We say a pair $(p, r)$ is admissible if $p, r \in [2,\infty]$ and $(p,r,d) \ne (2,\infty ,2)$ satisfying $\frac{2}{p}+\frac{d}{r}=\frac{d}{2}$ and $2\leq r\leq \frac{2d}{d-2}$.

\begin{lem}\label{str}
Let $(p,r)$ and $(\gamma,\rho)$ be two admissible pairs and let $\gamma',\rho'$  be conjugates of $\gamma$ and $\rho$. Then
\begin{itemize}
\item[(1)] for every $\phi\in L^2(\RR^d)$, the function $t\mapsto S_t\phi$ belongs to $L^p(\RR;L^r(\RR^d))\cap L^{\infty}(\RR;L^2(\RR^d))$ and there exists a constant $C$  such that
\begin{align}\label{eq-stri-1}
\|S_\cdot\phi\|_{L^p(\RR;L^r(\RR^d))}\leq C\|\phi\|_{L^2(\RR^d)}.
\end{align}
\item[(2)] Let $I$ be an interval of $\RR$ and $J=\bar{I}$ with $0\in J$. Then for every $f\in L^{\gamma'}(I;L^{\rho'}(\RR^d))$, the function $t\mapsto \Phi_f(t)=\int_0^tS_{t-s}f(s)\d s$ belongs to $L^p(I;L^r(\RR^d))\cap L^{\infty}(J;L^2(\RR^d))$ and there exists a constant $C$ independent of $I$  such that
\begin{align}
&\|\Phi_f\|_{L^{\infty}(J;L^2(\RR^d))}\leq C\|f\|_{L^{\gamma'}(I;L^{\rho'}(\RR^d))} ;\label{det-stri-ine-L2}\\
&\|\Phi_f\|_{L^p(I;L^r(\RR^d))}\leq C\|f\|_{L^{\gamma'}(I;L^{\rho'}(\RR^d))} .\label{det-stri-ine}
\end{align}
\item[(3)] For all $q\geq 2$ and all  $\mathbb{F}$-predictable process $\xi:[0,T] \times L^{2}(\RR^d)\times\Omega \rightarrow L^{2}(\RR^d;\mathbb{C})$ in \ $L^{q}\big(\Omega;L^{2}([0,T] \times L^{2}(\RR^d);L^{2}(\RR^d))\cap L^{q}([0,T] \times L^{2}(\RR^d);L^{2}(\RR^d))\big)$, we have
\begin{align}\label{str-ineq-stoch}
\EE\Big\|\int_0^\cdot\int_{L^2(\RR^d)} S_{\cdot-s}\xi(s,z)\tilde{N}(\d z,\d s)\Big\|^q_{L^p(0,T;L^r(\RR^d))}\leq &C_{q}\,\EE\Big(\int_0^T\int_{L^2(\RR^d)}\|\xi(s,z)\|^2_{L^2(\RR^d)}\nu(\d z)\d s\Big)^{\frac{q}{2}}\nonumber\\
&+C_{q}\,\EE\Big(\int_0^T\int_{L^2(\RR^d)}\|\xi(s,z)\|^q_{L^2(\RR^d)}\nu(\d z)\d s \Big).
\end{align}
\end{itemize}
\end{lem}

 %---------------------------
\section{ The proof of Theorem \ref{main}}
In this section, we shall prove the global existence of the original equation (\ref{NLS}). To do that, we first use the fixed point argument to show the global well-posedness of a truncated solution. Next we establish some uniform $L^2$-norm estimate of the local solution which combining with some stopping time arguments would lead us to believe that solutions ought to exist globally.

{\bf  Step 1.  Existence and uniqueness of a local solution:} First we define a truncation function $\theta$. Let $\theta:\mathbb{R}_+\rightarrow[0,1]$ be a non-increasing $C_0^{\infty}$ function such that
  $1_{[0,1]}\leq \theta \leq 1_{[0,2]}$ and $\inf_{x\in\RR_+}\theta'(x)\geq-2$. For $R\geq 1$, set $\theta_R(\cdot)=\theta(\frac{\cdot}{R})$.

Let us fix $R \geq 1$. Applying Lemma \ref{str}, similar to the proofs of Proposition 3.1, Proposition 3.6 and Proposition 3.7 in \cite{BLZ}, we can prove the existence and uniqueness of the global solution $Z^R$ to the following truncated equation
\begin{equation} \label{NLS-trucated}
     Z^R(t)=S_t x-i\lambda\int_0^tS_{t-s}\big(\theta_R(\|Z^R\|_{Y_s})|Z^R(s)|^{\al-1}Z^R(s) \big)\d s +\int_0^t\int_B S_{t-s}z\tilde{N}(\d z,\d s),\quad 0\leq t\leq T.
\end{equation}

%where $\gamma = \gamma\big(q,T, \|\hbar\|_{L^2(\RR^d)}, \int_{Z_1^c}\|z\|_{L^2(\RR^d)}^2 \nu(dz)\big)$ is independent of $R$.

We define a stopping time $\tau_R$ by
\[ \tau_R = \inf\{t \in [0,T]: \|Z^R\|_{Y_s}> R \}, \]
with the usual convention $\inf \emptyset =T$. Define $X(t)=Z^R(t)$ for $t\in[0,\tau_R]$ and $\tau_{\infty}=\lim_{R \uparrow\infty}\tau_R$. Then $(X(t))_{t\in[0,\tau_{\infty})}$ is a  maximal local mild solution of \eqref{NLS}.

{\bf Step 2. $L^2$-norm estimate:} For any $q \geq 2$, we need to prove 
\begin{equation}\label{L2 estimate}
\EE \sup_{t \in [0,T]}\|Z^R(t)\|_{L^2(\RR^d)}^q \leq \gamma,
\end{equation}
where $\gamma = \gamma\big(q,T, \|x\|_{L^2(\RR^d)}\big)$ is independent of $R$. Thus,  if $\PP(\tau_{\infty} = T) = 1$, then we can obtain 
(\ref{L2 norm}) by taking $R \rightarrow \infty$.

 As a weak equation in $H^{-1}(\RR^d)$, for $ 0\leq t\leq T$, we have, $\mathbb{P}$-a.s.
 \begin{equation*}
Z^R(t) = x +i \int_0^t \De Z^R(s) ds - i\lambda \int_0^t { \theta_R(\|Z^R\|_{Y_s})|Z^R(s)|^{\al-1}Z^R(s) ds} + \int_0^t \int_B z\tilde{N}(\d z,\d s).
\end{equation*}
Applying the It\^o formula, we get
\begin{eqnarray*}
&&\|Z^R(t)\|_{L^2(\RR^d)}^2 = \|x\|_{L^2(\RR^d)}^2  + \int_0^t{\int_B \|Z^R(s-) + z\|_{L^2(\RR^d)}^2 - \|Z^R(s-)\|_{L^2(\RR^d)}^2 \tilde{N}(dz ,ds)} \notag\\
&+& \int_0^t \int_B  \|Z^R(s) + z\|_{L^2(\RR^d)}^2 - \|Z^R(s)\|_{L^2(\RR^d)}^2 - 2 Re \langle Z^R(s), z\rangle_{L^2(\RR^d)} \nu(dz) ds,
\end{eqnarray*}
  where we used the fact that $Re\langle Z^R(s),i\De Z^R(s)\rangle_{L^2}=0$, since $i\Delta$ is skew-self-adjoint in $L^2(\RR^d)$.

Now we use the It\^{o} formula for the real-valued process $\Big(\|Z^R(t)\|^2_{L^2(\RR^d)}\Big)^{q/2}$ to get
\begin{eqnarray} \label{Ito}
&&\|Z^R(t)\|_{L^2(\RR^d)}^q = \|x\|_{L^2(\RR^d)}^q  + \int_0^t{\int_B \|Z^R(s-) + z\|_{L^2(\RR^d)}^q - \|Z^R(s-)\|_{L^2(\RR^d)}^q \tilde{N}(dz ,ds)} \notag\\
&+& \int_0^t \int_B  \|Z^R(s) + z\|_{L^2(\RR^d)}^q - \|Z^R(s)\|_{L^2(\RR^d)}^q - q \|Z^R(s)\|_{L^2(\RR^d)}^{q-2}Re \langle Z^R(s), z\rangle_{L^2(\RR^d)} \nu(dz) ds \notag \\
&+& \int_0^t \int_B \|Z^R(s)\|_{L^2(\RR^d)}^{q-2} \|z\|_{L^2(\RR^d)}^2 \nu(dz)ds.
\end{eqnarray}
By the elementary calculus, we have the following fact. For any $a, b \in L^2(\RR^d)$,
\begin{equation}\label{fact}
\big|\|a + b \|_{L^2(\RR^d)}^q - \|b\|_{L^2(\RR^d)}^q \big| \leq C_q (\|a\|_{L^2(\RR^d)}^q + \|b\|_{L^2(\RR^d)}^q). 
\end{equation}
Hence, by the Burkholder-Davis-Gundy inequality and Gronwall's inequality and (\ref{fact}), we obtain (\ref{L2 estimate}).

{\bf Step 3.  Global solution:} In order to obtain global solution of \eqref{NLS}, we only need to prove $\PP(\tau_{\infty} = T) =1$.
For the simplicity of presentation, we shall adopt the following notations for $t \in [0,T]$,
\begin{align*}
[\Psi^{R}(X)](t)=& -i\lambda\int_0^tS_{t-s}\big(\theta_R(\|Z^R\|_{Y_s})|Z^R(s)|^{\al-1}Z^R(s) \big)\d s, \\
M(t) =&  \int_0^t\int_B S_{t-s}z\tilde{N}(\d z,\d s).
\end{align*}
Applying (3) of Lemma \ref{str}, for all $q\geq2$,
\begin{eqnarray}\label{martingale}
\mathbb{E}(\|M\|^q_{L^p(0,T;L^r(\RR^d))})
&\leq& C_{q,T}\Big(\int_B \|z\|_{L^2(\RR^d)}^q \nu(d z)+\Big(\int_B \|z\|_{L^2(\RR^d)}^2 \nu(d z)\Big)^{q/2}\Big)\nonumber\\
&\leq& C_{q,T}\Big(\int_B\|z\|_{L^2(\RR^d)}^2 \nu(d z)+\Big(\int_B\|z\|_{L^2(\RR^d)}^2 \nu(d z)\Big)^{q/2}\Big) < \infty.
\end{eqnarray}
Let us fix $\om \in \Om$ and take $T_R(\om) \in (0, T]$ whose value will be determined later on.
By Lemma \ref{str} and \cite[Proposition 3.1]{BLZ}, we have
\ba
&&\|Z^R(\om)\|_{L^p(0,T_R;L^r(\RR^d))} \notag\\
&\leq& \|S_t x\|_{L^p(0,T_R;L^r(\RR^d)} +\|\Psi^{R}(Z^R(\om))\|_{L^p(0,T_R;L^r(\RR^d)} +\|M(\om)\|_{L^p(0,T_R;L^r(\RR^d)} \notag \\
&\leq&  C\|x\|_{L^2(\RR^d)} + CT_R^{1-\frac{(\al-1)d}{2}}\|Z^R(\om)\|_{L^p(0,T_R;L^r(\RR^d)}^{\al } + \|M(\om)\|_{L^p(0,T_R;L^r(\RR^d)} \notag\\
&\leq& M_R^{T_R}(\om)+ CT_R^{1-\frac{(\al-1)d}{4}}\|Z^R(\om)\|_{L^p(0,T_R;L^r(\RR^d)}^{\al },
\ea
where
\begin{equation}
 M_R^{T_R}(\om) := C\sup_{0\leq t < T_R} \|Z^R(t,\om)\|_{L^2(\RR^d)} + \|M(\om)\|_{L^p(0,T_R;L^r(\RR^d)}.
\end{equation}
Also let us denote
\begin{equation}
 M_R(\om) := C\sup_{0\leq t \leq T} \|Z^R(t,\om)\|_{L^2(\RR^d)} + \|M(\om)\|_{L^p(0,T;L^r(\RR^d)}.
\end{equation}
The following fact will be used. Fix $K>0$, consider a function $f$ defined by
$$
f(x)=K + \frac{x^{\al}}{4(2K)^{\al -1}} -x,\ x\geq0.
$$
Then there exist two points $c_1,c_2$ with $0<c_1<2M<c_2<4^{\frac{1}{\alpha-1}}\cdot 2K<\infty$ at which $f(c_1)=f(c_2)=0$, and
$f(x)\geq0$ if and only if $0\leq x\leq c_1$ or $x\geq c_2$. The proof of this fact can be found around (4.8) in \cite{BLZ}.
%the following fact
%\begin{equation}
%\forall x>0, \ \exists c_1\leq 2K,\ c_2 >c_1: \quad x \leq K + \frac{x^{\al}}{4(2K)^{\al -1}} \Rightarrow x \leq 2K \ or \ x> c_2,
%\end{equation}

Let $T_R(\om) = T \land (4C(2M_R(\om))^{\al-1})^{-\frac{1}{1-\frac{(\al-1)d}{4}}}$, then
\[
C T_R(\om)^{1-\frac{(\al-1)d}{4}}\big( C\sup_{0\leq t \leq T} \|Z^R(t,\om)\|_{L^2(\RR^d)} + \|M(\om)\|_{L^p(0,T;L^r(\RR^d)}\big)^{\al -1} \leq  \frac{1}{2^{\al +1}}. \]
Combining this with the fact above, we have
\begin{equation}\label{control}
\|Z^R(\om)\|_{L^p(0,T_R;L^r(\RR^d))} \leq2  M_R^{T_R}(\om)\leq 2 M_R(\om).
\end{equation}

This inequality suggests that we can replace $T_R$ which is not a stopping time by a stopping time $\si_R$ defined by
\begin{equation}
\si_R := \inf\{t \in [0, T]: C t^{1-\frac{(\al-1)d}{4}}\big( C\sup_{0\leq s < t} \|Z^R(s)\|_{L^2(\RR^d)} + \|M\|_{L^p(0,t;L^r(\RR^d)}\big)^{\al -1} > \frac{1}{2^{\al +1}} \}.
\end{equation}
Therefore, similar to the argument above, (\ref{control}) holds with $T_R$ replaced by $\si_R$.
Then we define a sequence $\si_R^j$ for $j =1,\cdots$, as follow. For $j=1$ we put $\si_R^1=\si_R$.
\begin{align}
\si_R^{j+1} := \inf\{t \in (\si_R^j, T]: 
C (t-\si_R^j)^{1-\frac{(\al-1)d}{4}}\big( C\sup_{\si_R^j \leq s < t} \|Z^R(s)\|_{L^2(\RR^d)} + \|M\|_{L^p(\si_R^j,t;L^r(\RR^d)}\big)^{\al -1} > \frac{1}{2^{\al +1}}\}.
\end{align}
By the definition of $\si_R^j$, we infer $\si_R^{j+1}-\si_R^j \geq T_R$, for each $j$. Let $N = \lfloor\frac{T}{T_R} \rfloor$. Then we can see that $\si_R^{N+1} = T$. Hence, applying similar arguments as above, for $j=1,\cdots, N$,
\[
\|Z^R\|_{L^p(\si_R^j,\si_R^{j+1};L^r(\RR^d))} \leq 2 M_R.
\]
Then,
\ba
\|Z^R\|_{L^p(0,T;L^r(\RR^d))} &\leq& \|Z^R\|_{L^p(0,\si_R^1;L^r(\RR^d))} + \sum_{j=1}^N \|Z^R\|_{L^p(\si_R^j,\si_R^{j+1};L^r(\RR^d))}
\leq 2( \frac{T}{T_R} +1)M_R \notag \\
&=& 2M_R +2T(4C(2M_R)^{\al-1})^{\frac{1}{1-\frac{(\al-1)d}{4}}}M_R \notag\\
&\leq& 2M_R + C_{T, \al,d}(M_R)^{\theta},
\ea
where $\theta = \frac{4(\al-1)}{4-(\al-1)d} +1$.
By applying (\ref{L2 estimate}) and (\ref{martingale}), we obtain
\ba \label{stnorm}
&&\EE\|Z^R\|_{L^p(0,T;L^r(\RR^d))}
\leq 2\EE M_R +C_{T, \al,d}\EE(M_R)^{\theta} \notag \\
&\leq& C_{T, \al,d,\theta}\big[\big(\EE\sup_{0\leq t \leq T} \|Z^R(t)\|_{L^2(\RR^d)} +\EE\sup_{0\leq t \leq T} \|Z^R(t)\|_{L^2(\RR^d)}^{\theta} \big) + 1\big]\notag \\
&\leq& C_{T, \al,d,\theta, \gamma}.
\ea

Combining with (\ref{L2 estimate}) and (\ref{stnorm}), we deduce that
 \begin{align*}
  \PP(\tau_R=T)&= \PP\{ \sup_{0\leq t\leq T}\|Z^R(t)\|_{L^2(\RR^d)} + \|Z^R\|_{L^p(0,T;L^r(\RR^d))} \leq R\}\\
  &\geq 1-\frac{\EE\sup_{0\leq t\leq T} \|Z^R(t)\|_{L^2(\RR^d)}+\EE \|Z^R\|_{L^p(0,T;L^r(\RR^d))}}{R}\geq 1- \frac{C_{T, \al,d,\theta, \gamma}}{R}.
 \end{align*}
Hence we have
\begin{align*}
 \PP(\tau_{\infty}=T)\geq \lim_{R \uparrow\infty}\RR(\tau_R=T)=1.
\end{align*}
Thus we infer that $\tau_\infty \geq T$,  $\PP$-a.s.. Since $T$ is arbitrary,
this shows $X(t)$, $t\in[0,\infty)$ is a unique global mild solution to (\ref{NLS}). %That is Assumption \ref{ass1} hold.

 %\bibliographystyle{elsarticle-num} 
 %\bibliography{cas-refs}


\end{document}